\documentclass{amsart}
\usepackage{graphicx}
\newtheorem{thm}{Theorem}[section]
\newtheorem{cor}[thm]{Corollary}
\newtheorem{lem}[thm]{Lemma}

\newtheorem{stat}[thm]{Statement}
\theoremstyle{definition}
\newtheorem{defn}[thm]{Definition}
\theoremstyle{remark}

\numberwithin{equation}{section}
\newcommand{\abs}[1]{\vert#1\vert}

\begin{document}

\title[]{Whitney Regularity Of The Image Of The Chevalley Mapping }%
\author{ G\'erard P. BARBAN\c CON}%
\address{University of Texas at Austin}%
\email{gbarbanson@yahoo.com }%

\thanks{I thank Professor Hubert Rubenthaler, University of Strasbourg, France who invited me to present this work in his seminar.}%

%\date{October 2009}%

% ----------------------------------------------------------------
\begin{abstract}
 A compact set $K\subset{\mathbb R}^n $ is Whitney 1-regular if 
 the geodesic distance in $K$ is equivalent to the Euclidean distance.
 Let $P$ be the Chevalley map defined by an integrity basis of
 the algebra of polynomials invariant by a reflection group, this note
 gives the Whitney regularity of the image by $P$ of any closed ball centered at the origin in ${\mathbb R}^n$.
 The proof uses the works of Givental', Kostov and Arnold on the symmetric group. It needs a generalization of a property of the Van der Monde determinants to the Jacobian of the Chevalley mappings.
\end{abstract}
\maketitle
% ----------------------------------------------------------------

\section{Introduction}

 \begin{defn}{\rm [13]}  A compact set $K \subset {\mathbb R}^n$,
  connected by rectifiable arcs, is Whitney 1-regular or has the Whitney
 property ${\mathbf P}_1$ if the geodesic distance in $K $  is
 equivalent to the Euclidean distance. \end{defn}

\noindent In other words, $K$ is 1-regular if there exists a constant $C_K >0 $, such
 that any couple of points $x$ and $x'$ in $K $ may be joined by a rectifiable arc in $K $ with
 length $\ell(x,x') \leq C_K | x-x'| $.
 \medskip

Let $W\in O(n)$ be a finite reflection group. The
algebra of $W$-invariant polynomials is generated by $n$ algebraically independent
$W$-invariant homogeneous polynomials [5] and
the degrees of these polynomials are uniquely determined [6]. Let
$p_1, \ldots , p_n$ be a set of such basic invariants and define the `Chevalley' mapping
$$P:{\mathbb R}^n\ni  x\mapsto (p_1(x), \ldots , p_n(x))\in {\mathbb R}^n.$$

This note gives a proof of the Whitney property ${\mathbf P}_1$ for the image by $P$ of any closed ball centered at the origin in ${\mathbb R}^n$.

For $I_2(p)$ a dihedral group, the 1-regularity of $P(\mathbb R^2)$ is clear since it is the salient part
cut in the plane by a cusp. In particular the property holds for the group $G_2=I_2(6)$.
Since there is nothing to prove for $I_2(p)$, it will not be considered in what follows.
The interested reader will find some details in appendix $A$.

For the symmetric group the property ${\mathbf P}_1$ of the image of the Newton mapping was
conjectured in [3] with a proof in dimension less than or equal to 4. The general proof was given
by Kostov [9] for the symmetric groups and for $A_n$ groups .

First we show that the results given for the symmetric group and announced for reflection groups by Givental in [7], actually hold for these groups.

Then we derive the property ${\mathbf P}_1$ for the image of the Chevalley map by using
the methods in [1] or [9], and already in [3].

This property is not altered by diffeomorphism, so it does not depend on
the choice of coordinates nor on the choice of the set of basic invariants
since a change of basic invariants is an invertible polynomial map on the target.
\vskip 5pt

The image sets of closed balls centered at the origin by Chevalley maps associated with reflection groups are natural but not trivial examples of sets with this geometric property.\vskip 5pt

Finally, the property ${\mathbf P}_1$ shares with the convexity an interesting feature. Let ${\mathcal O}$ be an open set of ${\mathbf R}^n$ such that its closure ${\mathcal F}=\overline {{\mathcal O}}$ has property ${\mathbf P_1}$ and let $f\in C^m({\mathcal O})$. If the field of Taylor polynomials of $f$ has a continuous extension of degree $m$ on ${\mathcal F}$, this extension is Whitney $m$-regular and by the Whitney extension theorem [12] there is one and actually there are infinitely many extensions of $f$ on ${\mathbb R}^n$.

\begin{thm} {\rm [13]}  Let ${\mathcal O}$ be an open set of ${\mathbf R}^n$ and assume that its closure ${\mathcal F}=\overline {{\mathcal O}}$ has property ${\mathbf P_1}$. If $f\in C^m({\mathcal O})$is such that whenever $\mid k\mid=m$, the $\partial^{\mid k\mid}f/ \partial x^k$ have a continuous extension to ${\mathcal F}$, then $f$ has an extension $\tilde f\in C^m({\mathbb R}^n)$. \end{thm}

 Thanks to the 1-regularity of $P({\mathbb R^n})$, the proof of the theorem of Chevalley in finite class of differentiability [4], would be easier. Let $f\in C^r({\mathbb R}^n)^W$ be $W$ invariant and $f=F\circ P$. The study of the continuity of the derivatives of $F$ on the border of $P({\mathbb R^n})$, shows that the field of Taylor polynomials of $F$ has a continuous extension of degree $[r/h]$ on $P({\mathbb R^n})$, where $h$ is the Coxeter number of $W$. By theorem 1.2 there is one, and actually there are infinitely many extensions of $F$, of class $C^{[r/h]}$ on ${\mathbb R}^n$.

\section{The Chevalley mapping, strata of its image and minors of its Jacobian}

  When $W$ is reducible, it is a direct product of its irreducible
 components, say $W= W^0\times W^1 \times \ldots \times W^s $ and we may
 write $\mathbb{R}^n$ as an orthogonal direct sum $\mathbb{R}^{n_0}\oplus
 \mathbb{R}^{n_1}\oplus \ldots \oplus \mathbb{R}^{n_s}$. The first component $W^0$ is
 the identity on the $W$-invariant subspace $\mathbb{R}^{n_0}$,
 and for $i=1,\ldots, s$, $W^i$ is a finite Coxeter group
 acting on $\mathbb {R}^{n_i}$.

 For an irreducible $W$ we will assume as we may that the degrees of the
 polynomials $p_1,\ldots ,p_n$ are in increasing order:
 $2=k_1\le \ldots \le k_n=h$, where $h$ is the Coxeter number of $W$, even more that $p_1(x)=\sum_{i=1}^n x_i^2 =|x|^2$.
 
 In the reducible case, using coordinates that fit with the orthogonal direct sum,
 if $w=w_0 w_1\ldots w_s\in W$ with $w_i\in W^i, \; 0\le i\le s$,
 for all $x\in \mathbb{R}^{n}$ we have
 $w(x)= w(x^0, x^1,\ldots , x^s)= (x^0, w_1(x^1),\ldots, w_s(x^s))$ . The
 direct product of the identity $P^0$ on $\mathbb {R}^{n_0}$ and Chevalley
 mappings $P^i$ associated with $W^i$ acting on $\mathbb{R}^{n_i},\; 1
 \le i\le s$, is a Chevalley map associated with the action of $W$
 on $\mathbb{R}^{n}$.
 For each $i=1\ldots s$ the $p^i_j$ are of degree $k^i_j$, for $j=1$ to $n_i$.
 The $k^i_1$ are all equal to 2, and $k^i_{n_i}$ is the Coxeter number of $W^i$.

 If for each $i$, $P^i(\mathbb{R}^{n_i})$ is 1-regular, so is
 $P(\mathbb{R}^{n})$. As a consequence it is sufficient to prove the regularity
 when $W$ is irreducible and from now on, we will assume $W$ to be
 a finite Coxeter group.

Let $\mathcal{R}$ be the set of reflections different from identity
in $W$. The number of these reflections is
$\mathcal{R}^{\#}=d=\sum_{i=1}^n(k_i-1)$. For each $\tau \in
\mathcal{R}$, we consider a linear form $ \lambda_{\tau}$ the
kernel of which is the reflecting hyperplane
$ H_{\tau}=\{ x\in {\bf R}^n |\tau (x) = x \} $.

The jacobian determinant of $P$ is $\abs{J_P}=c
\prod_{\tau \in \mathcal{R}} \lambda_{\tau}$ for some constant
$c\neq 0$. The critical set is the union of the $ H_{\tau}$ when
$\tau$ runs through $\mathcal{R}$.

 A Weyl Chamber $C$ is a connected component of the regular
 set. The other connected components are obtained by the
 action of $W$ and the regular set is $\bigcup_{w\in W} w(C)$.
 There is a stratification of ${\mathbb R}^n$ by the regular set, the
 reflecting hyperplanes $ H_{\tau}$ and their intersections.
 The mapping $P$ is neither injective
 nor surjective but it induces an analytic diffeomorphism of $C$ onto
 the interior of $P({\mathbb R}^n)$ and an homeomorphism
 that carries the stratification from the
 fundamental domain $\overline{C}$ onto $P(\mathbb{R}^n)$.

 In the fundamental domain $\overline{C}$, the walls of the Weyl chamber $C$ are contained in $n$ hyperplanes $(H_{\omega})_{\omega \in \Omega}$, where $\Omega$ is a subset of cardinal $n$ in $\mathcal {R}$.

 A stratum $S$ of dimension $k$ in $\overline{C}$ is the intersection of $(n-k)$ of the $H_{\omega}$.
 The $\lambda_{\tau}$ that are linear combinations of the $\lambda_{\omega}$ vanishing on $S$ will also vanish there, so that $p \ge n-k$ hyperplanes $H_{\tau}$ will intersect along $S$.

 The points in $S$ have the same isotropy subgroup $W_{S}$, generated by the reflections in the $p$ hyperplanes $H_{\tau}$ containing $S$. In a neighborhood of $S$, since $P$ is $W_S$ invariant we can write $P=Q\circ V$, where $Q$ is invertible and $V$ is a Chevalley mapping for $W_S$.

 \noindent We write $W_S= W^0\times W^1 \times \cdots \times W^\ell $ where $W^0$ is the identity on $S$ and the $W^m, 1\le m\le \ell$ are the other irreducible components of $W_S$.
 We choose coordinates fitted to the orthogonal direct sum
 $\mathbb {R}^n=\mathbb {R}^{k}\oplus \mathbb {R}^{n_1}\oplus \cdots \oplus \mathbb {R}^{n_\ell}$, and we have
 $V={\rm Id}_k\times V_1\times \cdots\times V_\ell$, where $V_m$ is a Chevalley map associated with $W^m$.

\noindent The jacobian matrix $J_V$ is block diagonal. Let $I_k,
J_{V_1}, \ldots, J_{V_\ell}$ be the diagonal blocks, the
determinants $\abs{J_{V_1}}, \ldots, \abs{J_{V_\ell}}$, all vanish
on $S$. Moreover, the rank of the $V_i$ on $S$ is $0$. The $k\times k$ minor $\abs{I_k}$ in the upper left corner
is $1$ but all the bordering $(k+1)\times (k+1)$ minors vanish on
$S$. The rank of $V$ on $S$, which is the rank of the Chevalley mapping $P$
on $S$, is $k$ the dimension of $S$.

 When restricted to $S$, $P_{|S}= Q\circ V_{|S}$ is an analytic isomorphism on its image
 and onto each of its projections on the spaces $\mathbb{R}^{\ell}, \; k\le \ell\le n-1$.

\section{The varieties $P_k^{-1}[P_k(x)]$.}

   For $1\le k\le n$ we set $P_k=(p_1, p_2, \ldots,p_k):\mathbb{R}^n\rightarrow \mathbb{R}^k$,
  and analogously $Q_k=(q_1, q_2, \ldots,q_k)$,  when $k=n$ however we may write $P$ as usual. We denote with $\Pi_k$ the image $P_k(\mathbb{R}^n)$, and by $\Pi_k'$ the image by $P_k$ of the ball centered at the origin of radius $a$, that is $\Pi_k\cap \{p|0\le p_1\le a^2\}$.\medskip

  Let $m^k\in \Pi_k$, then \[P_k^{-1}(m^k)=\{x\in \mathbb{R}^n; p_1(x)=m_1,\ldots, p_k(x)=m_k\}\] is a compact algebraic variety of co-dimension at most $k$.

  \begin{thm}\label{thm1} For almost all $m^k$, the intersection of a fundamental domain
  $\overline{C}$ with $P_k^{-1}(m^k)$ is either empty or contractible. \end{thm}

 By Sard's lemma, for almost all $m^k$, the variety $P_k^{-1}(m^k)$ is a non singular manifold of dimension $n-k$, transversal to all strata in $\overline{C}$.
 If $m^k$ is on the border of $\Pi_k$, the variety is singular. Even if $m^k$ belongs to
 the interior of $\Pi_k$, it may be the image of both regular and singular points. In this case it is
 a singular image and  $P_k^{-1}(m^k)$ is a singular manifold.

 \begin{lem}\label{lem1} For almost all $m_k$, $P_k^{-1}(m^k)$ is a non singular manifold and on this manifold, $p_{k+1}$ is a Morse function with no critical point on the strata of positive dimension in $P_k^{-1}(m^k)\cap\bar{C}$. Its critical points are the zero dimensional strata of this intersection, that is the non empty intersections of the manifold and the strata of dimension $k$ in $\overline{C}$. \end{lem}

\noindent See Lemma 1' - given in [7] without a complete proof.

 {\it Proof.} On the non critical manifold $P_k^{-1}(m^k)$, there is a $k\times k$ minor of the matrix $\partial (p_1,\ldots,p_k)/\partial (x_1,\cdots, x_n)$ that does not vanish and we may assume that it is the principal minor $\partial (p_1,\ldots,p_k)/\partial (x_1,\ldots, x_k)$ which is not 0.

 Let $L=p_{k+1}-\mu_k (p_k-m_k)-\ldots-\mu_1 (p_1 -m_1)$ be the Lagrange function.
The critical points of $p_{k+1}$ on $P_k^{-1}(m^k)$ are the points where:
\[{\partial p_{k+1}\over \partial x_i}-\mu_k {\partial p_{k}\over \partial x_i}-\ldots -\mu_1{\partial p_{1} \over \partial x_i}=0,\quad {\rm for}\; i=1,\ldots,n. \]
For this system to have solutions in the $\mu_i$s, $1\le i\le k$,
the $(k+1)\times (k+1)$ bordering minors
$\partial (p_1, \ldots, p_k,p_{k+1})/\partial(x_1, \ldots, x_k,x_{k+s}), \; s=1,\ldots,n-k$ must all be zero.

\begin{stat} \textit{If a $k\times k$ minor \[\partial (p_1, \ldots, p_k)/\partial(x_{i_1}, \ldots, x_{i_k})\] of the Jacobian of $P$ does not vanish but the $(k+1)\times (k+1)$ minors \[\partial (p_1, \ldots, p_k,p_{k+1})/\partial(x_{i_1}, \ldots, x_{i_k},x_{i_j})\] all vanish, then the rank of $P$ is $k$.} \end{stat}

 The statement, of interest by itself, is a consequence of the usual properties of the Van der Monde determinants for $A_n$ and $B_n$.

 It also holds for $D_n$ by the anti-symmetry of the minor determinants. We consider as we may the fundamental domain defined by the set of simple roots:
\[\alpha_n=e_{n-1}+e_n, \alpha_i=e_i-e_{i-1}, 1\le i\le n-1.\]
In this fundamental domain, $x_1\ge x_2\ge ...\ge x_{n-1}\ge |x_n|$, so if $x_1x_2...x_k=0$ then $x_{k+1} ... x_n$ is also 0. When $n$ is even, $n=2k$, we take as $k^{th}$ row in the jacobian determinant the row of derivatives of the $x_i^{n}$ and as $k+1^{st}$ row, the row of derivatives of $x_1x_2 \ldots x_n$.

For $H_3, H_4, D6$ and $F_4$, the result has been checked by computations with Mathematica \footnote{I am indebted to he University of Texas at Austin and Patrick Goetz from the Mathematics Department who provided me with the possibility of using Mathematica.}.

\textit{ Proof of the Statement}.

 First observe that if all the $(k+1)\times (k+1)$ minors  \[\partial
 (p_1,\ldots, p_{k+1})/\partial(x_{i_1}, \ldots, x_{i_k},x_{i_j})
 \] are $0$ for all $x_j$ not in $\{x_{i_1}, \ldots, x_{i_k}\}$, then all the $(k+1)\times (k+1)$ minors  \[\partial (p_1,\ldots, p_{k+1})/\partial(x_{j_1}, \ldots ,x_{j_{k+1}})\]also vanish.

 By Lagrange expansion  along the $k+2^{nd}$ row, the  \[\partial (p_1,\ldots, p_{k+2})/\partial(x_{i_1}, \ldots, x_{i_j},x_{i_l})\] vanish and so do all the $\partial (p_1,\ldots, p_{k+2})/\partial(x_{j_1}, \ldots, x_{j_{k+2}})$.

 By a finite induction, we get that the jacobian determinant itself is $0$ and the rank is at most $n-1$. So we are on the closure of a stratum of dimension $n-1$, which is a fundamental domain for a group, say
 $W_{|H_i}$, acting on the hyperplane $H_i$ of equation $\lambda_i(x)=0$ and generated by the reflections in the intersections of $H_i$ and the other $H_{\omega}$.

 We consider a set of simple roots orthogonal to the $H_\omega$. Let $y=(y_1,\ldots ,y_n)$ be the coordinates in this basis, we set $y=M(x)$ where $M$ is the matrix of change of basis, the rows of which are the $(\lambda_s)_{s=1}^{n}$. The matrix $M$ is known for any Coxeter group (see for instance ~[8]). We denote by $\overline{S_{n-1}^i}$ the closed stratum in $H_i$.

 Let $U=(u_1,\ldots,u_{n-1})$ be a system of basic invariants for $W_{|H_i}$.  Since $P$ restricted to $H_i$ is also invariant by this group, on $\overline{S_{n-1}^i}$ we have: $P_{|\overline{S_{n-1}^i}}=R\circ U\circ M$, where $R$ and $M$ are invertible.

  Denoting with $\hat{x}_i, \hat{y}_i$ the variables that do not show up, from the induction process we have: \[\partial (p_1,\ldots, p_{n-1})/\partial(x_1, \ldots,\hat{x_i},\ldots ,x_{n})=\] \[ \partial (p_1,\ldots, p_{n-1})/\partial(u_{1}, \ldots ,u_{n-1})\textbf{.} \partial (u_1,\ldots, u_{n-1})/\partial(y_{1}, \ldots,\hat{y_i},\ldots ,y_{n})\] \[\textbf{.}\partial (y_1,\ldots, \hat{y_i},\ldots, y_{n})/\partial(x_{1}, \ldots, \hat{x_i},\ldots,x_{n})=0.\] On one hand $\partial (p_1,\ldots, p_{n-1})/\partial(u_{1}, \ldots ,u_{n-1})\neq 0$ and on the other hand, for each Coxeter group there is an ordering of the simple roots, that is of the rows of the matrix $M$ such that for $i=1,\ldots, n$, $\partial (y_1,\ldots, \hat{y_i},\ldots, y_{n})/\partial(x_{1}, \ldots, \hat{x_i},\ldots,x_{n})\neq 0$, we then get: \[\partial (u_1,\ldots, u_{n-1})/\partial(y_{1}, \ldots,\hat{y_i},\ldots ,y_{n})=0\].

 As a consequence, for all $j=1,\ldots, n$: \[\partial (p_1,\ldots, p_{n-1})/\partial(x_1, \ldots \hat{x_j},\ldots ,x_{n})= \] \[\partial (p_1,\ldots, p_{n-1})/\partial(u_{1}, \ldots ,u_{n-1})\textbf{.}\partial (u_1,\ldots, u_{n-1})/\partial(y_{1}, \ldots,\hat{y_i},\ldots ,y_{n-1})\]\[\textbf{.}\partial (y_1,\ldots,\hat{y_i},\ldots ,y_{n-1}) /\partial(x_1, \ldots \hat{x_j},\ldots ,x_{n})=0,\] and  \[\partial (p_1,\ldots, p_{n-2},p_n)/\partial(x_1, \ldots \hat{x_j},\ldots ,x_{n})= \] \[\partial (p_1,\ldots, p_{n-2},p_n)/\partial(u_{1}, \ldots ,u_{n-1})\textbf{.}\partial (u_1,\ldots, u_{n-1})/\partial(y_{1}, \ldots,\hat{y_i},\ldots ,y_{n-1})\]\[\textbf{.}\partial (y_1,\ldots,\hat{y_i},\ldots ,y_{n-1}) /\partial(x_1, \ldots \hat{x_j},\ldots ,x_{n})=0.\]
Hence the rank is at most $n-2$ on $\overline{S_{n-1}^i}$ and just the same on all the $\overline{S_{n-1}^j}$ in the walls $H_j$ of $\overline{C}$, with $j=1,\ldots ,n$. So  we are on the closure of a stratum of dimension $n-2$, say $\overline{S_{n-2}^m}$.

The process may be iterated until we reach a stratum of dimension $k$ where the $k\times k$ determinant  \[\partial (p_1, \ldots, p_k)/\partial(x_{i_1}, \ldots, x_{i_k})\] is not $0$ by assumption. On this stratum the rank cannot be more nor less than $k$, it is $k$ and this achieves the proof of the statement.   $\square$

From Statement 3.3, the rank of $P$ is $k$ at the critical points of $p_{k+1}$ on $P_k^{-1}(m^k)$, these points
belong to strata $S$ of dimension $k$ in $\overline{C}$. They are points of 0-dimensional strata in the intersection $P_k^{-1}(m^k)\cap \overline{C}$.

If $P_k^{-1}(m^k)$ is not singular, strata of lower dimension do not intersect it since on such a stratum all $k\times k$ minors vanish. On strata of dimension $>k$, one of the $(k+1)\times (k+1)$ minors does not vanish and it must be one of the $$\partial (p_1, \ldots, p_k,p_{k+1})/\partial(x_1, \ldots, x_k,x_{k+s}), \; s=1,\ldots,n-k$$ since if not, by statement 3.3, the rank would be $\le k$. So $p_{k+1}$ may not have a critical point on those strata.

Let us show that the critical points of $p_{k+1}$ on $P_k^{-1}(m^k)$ are non degenerate.

\noindent By an orthogonal change of coordinates, we may assume that
$\{x_1,\ldots ,x_k\}$ are the coordinates in $S$ while $\{x_{k+1},\ldots ,x_n\}$ are the coordinates in $S^{\perp}$.

\noindent At a critical point, for $1\le i\le k, v_i=x_i$ and we have:
\[{\partial p_{k+1}\over \partial x_i}-\mu_k {\partial p_{k}\over \partial x_i}-\ldots
-\mu_1{\partial p_{1}\over \partial x_i}={\partial q_{k+1}\over \partial v_i}
-\mu_k {\partial q_{k}\over \partial v_i}-\ldots -\mu_1{\partial q_{1}\over \partial v_i} =0.\]
For $i\ge k$, on each $\mathbb{R}^{n_m},\; V^m$ is a Chevalley map for the irreducible component $W_S^m$ of $W_S$, and we have:
 \[{\partial p_{k+1}\over \partial x_i}-\mu_k {\partial p_{k}\over \partial x_i}-\ldots
-\mu_1{\partial p_{1}\over \partial x_i}=\sum_1^{n_m}\Big({\partial q_{k+1}\over \partial v_j}
-\mu_k {\partial q_{k}\over \partial v_j}-\ldots -\mu_1{\partial q_{1}\over \partial v_j}\Big)\Big({\partial v_j\over\partial x_i}\Big)
=0,\]
 but now ${\partial v_j/\partial x_i}=0 $ on $S$ since $v_j$ is an homogeneous polynomial of degree at least 2 in the $x_i$ for an $x\in \mathbb{R}^{n_m}$,
 and these $x_i$s vanish on $S$.

The vectors $(\partial P/\partial x_1,\ldots,\partial P/\partial x_k)$ span the $k$-dimensional plane tangent to $P(S)$.
The $\partial P/\partial x_s, k+1\le s\le n$ are linear combinations of them, hence the vanishing of the bordering minors.

The linearly independent vectors $\partial Q/\partial v_i$, $k+1\le i\le n$, span the orthogonal complement of the tangent plane. Accordingly the minors of order $\ge k+1$ in the jacobian of $Q$ do not vanish, and  for $i\ge k+1$,
\[{\partial q_{k+1}\over \partial v_i}-\mu_k {\partial q_{k}\over \partial v_i}-\ldots -\mu_1{\partial q_{1}\over \partial v_i}\neq 0.\]

When restricted to $S$, the second derivatives of $L$
\[{\partial^2 p_{k+1}\over \partial x_i\partial x_j}-\mu_k {\partial^2 p_{k}\over \partial x_i\partial x_j}-\ldots -\mu_1{\partial^2 p_{1}\over \partial x_i\partial x_j}\]
gets no contribution for the variables $x_i, x_j$ when $1\le i,j\le k$. For each $\mathbb{R}^{n_m}$, we get:
 \[\sum_{1\le r,t\le n_m}\Big({\partial^2 q_{k+1}\over
\partial v_{r}\partial v_t
 }-\mu_k {\partial^2 q_{k}\over \partial v_r\partial v_t}-\ldots -\mu_1
{\partial^2 q_{1}\over \partial v_r\partial v_t}\Big)\Big({\partial v_r \over\partial x_i}\Big)\Big({\partial v_t \over\partial x_j}\Big)+ \]
 \[\sum_{1\le r\le n_m}\Big({\partial q_{k+1}\over \partial v_r}-\mu_k {\partial q_{k}\over \partial v_r}-\ldots -\mu_1 {\partial q_{1}\over \partial v_r}\Big)\Big({\partial^2v_r\over \partial x_i\partial x_j}\Big).\]

 As already mentioned the first derivatives ${\partial v_r /\partial x_i}$ all vanish on $S$ and so do the mixed derivatives ${\partial^2v_r /{ \partial x_i\partial x_j}}$ either because $v_r$ does not depend on $x_i$ or $x_j$, or because it is an homogeneous polynomial of degree $> 2$ or a sum of squares of $x_i$s. For the first two of these reasons, many terms in the pure derivatives also vanish on $S$.
However for each $V_m:\mathbb{R}^{n_m}\rightarrow \mathbb{R}^{n_m}$, let $v_m^1$ be the first $W_S^m$-invariant which is the sum of the squares of the $x_i$ in $\mathbb{R}^{n_m}$. We have ${\partial^2v_m^1/ \partial x_i^2}=2$, so that, with the above remark:

\[{\partial^2 p_{k+1}\over \partial x_i^2}-\mu_k {\partial^2 p_{k}\over \partial x_i^2}-\ldots -\mu_1
{\partial^2 p_{1}\over \partial x_i^2}= 2({\partial q_{k+1}\over \partial v_m^1}-\mu_k
{\partial q_{k}\over \partial v_m^1}-\ldots -\mu_1{\partial q_{1}\over \partial v_m^1})\neq 0.\]

\noindent For each $W_S^m$, the quadratic differential is definite with the sign of \[{\partial q_{k+1}/ \partial v_m^1}-\mu_k {\partial q_{k}/ \partial v_m^1}-\ldots -\mu_1{\partial q_1/ \partial v_m^1}\] which may be different for different irreducible components. Overall the quadratic differential restricted to $S$, kernel of the first differential,
is of the form $Q=\bigoplus_{i=1}^l Q_i$ where $Q_i$ is a definite quadratic form on $\mathbb{R}^{n_i}$.\vskip 5pt

 Accordingly $p_{k+1}=q_{k+1}\circ V$ is a Morse function on $P_k^{-1}(m^k)$. $\square$

\noindent By the equivariant Morse lemma [2], in the neighborhood of a critical point at the intersection of
$S$ and $P_k^{-1}(m^k)$, the restriction of $p_{k+1}$ is $W_S$-locally equivalent to a $W_S$-invariant quadratic form which is the direct sum of definite quadratic forms in each $\mathbb{R}^{n_i}$. The pair $(S, p_{k+1})$ is locally equivalent to $(\mathbb{R}_+^u\oplus \mathbb{R}_+^v, Q_+\oplus Q_-) $.

 \begin{lem}\label{lem2}{\rm [7
 ]} The reconstruction of the topology of a level set of a function on the sum of quadrants $\mathbb{R}_{+}^a\oplus \mathbb{R}_{+}^b$ in the neighborhood of the critical point $0\oplus 0$ with the quadratic differential $Q_{+}\oplus Q_{-}$ is trivial if $a, b>0$ and consists of the birth (death) of a simplex otherwise.\end{lem}

\noindent Theorem ~\ref{thm1} may now be proved by induction on $k$, see [7].\medskip%

In particular, and this is what we need, for almost all $m^k\in \Pi_k,\; P_k^{-1}(m^k)\cap \overline{C}$ is connected.

\begin{cor}\label{cor1} Every variety $P_k^{-1}(m^k) \cap \overline{C},\; k=1,\ldots, n$, is connected or empty. \end{cor}

 \noindent The corollary may be derived from ~\ref{thm1}, exactly as in ~{[9]}. Basically if $m^k$ is not in the regular image and $P_k^{-1}(m^k)$ is of dimension $n-s$, on the intersection with a stratum of dimension $l$ the rank of $P_k$ is $s$ if $l\ge s$ and $l$ if $l\le s$. The intersection is of dimension $l-s$ in the first case and 0-dimensional in the second case. So, either they intersect transversally or the intersection is reduced to a finite number of points. When the intersection is transversal $P_k^{-1}(m^k)$ is connected since almost all the $P_k^{-1}(n^k)$ are, for $n^k$ close enough to $m^k$. The second case is taken care of by an induction on $k$, using the following:

 \begin{lem}\label{lem3}{\rm ~{[9]}} Let $K$ be a connected compact set in $\mathbb {R}^n$, and $f$ be a function continuous on $K$. If all but a finite number of the level sets of $f$ in $K$ are connected, then they are all connected.\end{lem}

 We observe that if not empty$, P_1^{-1}(m^1)$ intersection of $\overline{C}$ and a sphere, is connected and we assume that every non empty $P_k^{-1}(m^k) \cap \overline{C}$ is connected. If we study a $P_{k+1}^{-1}(m^{k+1})$ not covered by the previous cases, we apply ~\ref{lem3} to the compact connected $K=P_k^{-1}(m^k) \cap \overline{C}$ containing it and $f=p_{k+1}$ to reach the conclusion for the level set $P_{k+1}^{-1}(m^{k+1})$.  $\square$

  As a consequence, $P(P_k^{-1}(m^k))$ is connected or empty and we have:

\begin{cor}\label{cor2} The fibres of the projections: $\Pi_{k+1}\rightarrow \Pi_k$ are connected. They are points or intervals.\end{cor}

\section{Whitney 1-regularity of $P(\mathbb {R}^n)$}

 Let us consider the lift-up $(p_1,\ldots, p_k)\mapsto p_{k+1}$ of the projection $P_k(S)\subset \Pi_k$ of a stratum of dimension $k$.
 The derivatives are obtained by solving the system:
  \[ {\partial p_{k+1}\over \partial x_i}=\sum_{j=1}^k {\partial p_{k+1}\over \partial p_j}\;{\partial p_j\over \partial x_i},
  \quad i=1,\ldots, k\]
  by Cramer's method. Considering the subgroup of $W$ generated by the reflections in the hyperplanes that do not contain $S$,
  for $j=1,\ldots, k $ the $\partial p_{k+1}/ \partial p_j$, are quotients of two polynomials anti-invariant by this subgroup.
  So, they are rational fractions the numerator and the denominator of which are invariant homogeneous polynomials that do not
  vanish but at the origin. Since the degree of the numerator is greater than the degree of the denominator the rational fractions
  have continuous extensions on $\bar{S}$.

 \noindent If we restrict ourselves to the compact subset $\Pi_k'\subset P(\mathbb{R}^n)$, determined by
 $0\le p_1\le a^2$, the $x_i$s and as a consequence the $\partial p_{k+1}/ \partial p_j,\; j=1,\ldots, k $ are
 bounded on $\overline {S}\cap \Pi_k'$.

 \noindent $P$ is a homeomorphism of $\bar{S}$ onto its image $P(\bar{S})$, and  so is $P_k$ (on any compact it
 is continuous and one to one). Hence $p_{k+1}$ which is continuous with respect to the variables $(x_1,\ldots,x_k)$
 on $\bar{S}$, is also continuous in the variables $(p_1,\ldots,p_k)$ on $P_k(\bar{S})$ and moreover by the previous paragraph it is Lipschitz on $\overline {S}\cap \Pi_k'$.

 The border of $\Pi_{k+1}$ is contained in the images $P_{k+1}(\bar{S})$ of closed strata of dimension $k$.
 These images are graphs of functions $p_{k+1}$ on $P_k(\bar{S})$. By ~\ref{cor2}, the graph of one of the $p_{k+1}$,
 say $p_{k+1}^{\rm max}$, is above the others and another, say the graph of $p_{k+1}^{\rm min}$, is below the others.
 In $\Pi_k$, the images of the closure of strata of dimension $k$ will intersect along the images of strata of lesser
 dimension. Above such points of intersection the mapping $p_{k+1}^{\rm max}$ (resp. $p_{k+1}^{\rm min}$) may and will
 change but globally $p_{k+1}^{\rm max}$ (resp. $p_{k+1}^{\rm min}$) will still be continuous and Lipschitz since the
 functions above the two strata are glued by their common value above the stratum of lower dimension along which they intersect.

Now, it is easy to get for any two points $A$ and $B$ in the compact $K= P(\mathbb{R}^n)\cap \{p|0\le p_1\le a^2\}$ a continuous arc
$AB\subset K$ of length $\ell(AB) \leq k_K \abs{AB}$, following the method in ~{[9]} and already in ~{[3]}.

We may also observe with ~{[1]} that the prism between graphs of Lipschitz functions over a Whitney
1-regular domain is Whitney 1-regular. Since
$\Pi_1'$ is Whitney 1-regular, by induction, assuming $\Pi_{k}'$ to be Whitney 1-regular, the prism $\Pi_{k+1}'$ over it and between the graphs of the Lipschitz functions $p_{k+1}^{\rm min}$ and $p_{k+1}^{\rm max}$ will also be Whitney 1-regular.
Hence all the $\Pi_k'$ and in particular $\Pi_n'=P(\mathbb{R}^n)\cap \{p|0\le p_1\le a^2\}$ are Whitney 1-regular and we can state:

\begin{thm} The image by the Chevalley mapping of any closed centered at the origin of $\mathbb{R}^n$ is Whitney 1-regular.

\end{thm}

% ----------------------------------------------------------------

% ----------------------------------------------------------------

\end{document}